\documentclass[11pt, twoside, eqno]{article}
\usepackage{latexsym}
\usepackage{amssymb}
\usepackage{amsfonts}
\begin{document}
\begin{center}

\noindent {\bf \Large On the paper "On Green's relations, 
$2^0$-regularity and Quasi-ideals in
$\Gamma$-semigroups"}\bigskip

\medskip

{\bf Niovi Kehayopulu}\bigskip

{\small \it Department of Mathematics,
University of Athens \\
157 84 Panepistimiopolis, Athens, Greece \\
email: nkehayop@math.uoa.gr}\medskip

September 9, 2013 \end{center}\bigskip

\bigskip

\noindent This is about the paper in the title published in by
Hila and Dine in [1]. According to Lemmas 2.13 and 2.14 of this
paper a $\Gamma$-semigroup $M$ is left (resp. right) regular if
and only if $a{\cal L} (a\gamma a)$ (resp. $a{\cal R} (a\gamma
a)$) for all $a\in M$ and all $\gamma\in\Gamma$. The Theorem 2.15,
Corollary 2.16, Proposition 2.17, Proposition 2.21, Theorem 2.22,
Corollary 2.23, Theorem 2.24, Theorem 2.25, Proposition 3.7,
Proposition 4.3, Lemma 4.4, Theorem 4.5, and the Corollary 4.7 of
the paper are based on these Lemmas. For Lemmas 2.13 and 2.14 by
Hila and Dine in [1], the authors refer to Theorem 1 in the paper
by Kwon and Lee in [2], but the proof of Theorem 1 in [2] is
wrong.

The authors use the usual definition of left (right) regular
$\Gamma$-semigroups. According to Definition 2.11 of this paper a
$\Gamma$-semigroup $M$ is called {\it left regular} if $a\in
M\Gamma a\Gamma a$ for every $a\in M$. That is, for every $a\in M$
there exist $x\in M$ and $\gamma, \mu\in\Gamma$ such that
$a=x\gamma a\mu a$. According to Definition 2.12, a
$\Gamma$-semigroup $M$ is called {\it right regular} if $a\in
a\Gamma a\Gamma M$ for every $a\in M$. That is, for every $a\in M$
there exist $x\in M$ and $\gamma, \mu\in\Gamma$ such that
$a=a\gamma a\mu x$. In [2] the authors consider the following
definition of a $\Gamma$-semigroup (shortly
$po$-$\Gamma$-semigroup): Let $M$ and $\Gamma$ be two nonempty
sets. $M$ is called a {\it $\Gamma$-semigroup} if

(1) $M\Gamma M\subseteq M$, $\Gamma M\Gamma\subseteq \Gamma$ and

(2) $(a\gamma b)\mu c=a(\gamma b\mu)c=a\gamma (b\mu c)$ for all
$a,b,c\in M$ and all $\gamma,\mu\in\Gamma$.\\
Then an ordered $\Gamma$-semigroup is a $\Gamma$-semigroup endowed
with an
order relation "$\le$".\\
This is the Theorem 1 of the paper in [2]:

\noindent {\bf Theorem 1.} {\it Let M be a
$po$-$\Gamma$-semigroup. The following are equivalent:

$(1)$ M is left regular.

$(2)$ $L(a)\subseteq L(a\gamma a)$ for every $a\in M$ and every
$\gamma\in\Gamma$.

$(3)$ $a{\cal L} (a\gamma a)$ for every $a\in M$ and every
$\gamma\in\Gamma$.}\\
And this is the proof of $(1)\Rightarrow 2$ of Theorem 1 (I copy
it from [2]): "Let $M$ be left regular. If $t\in L(a)$, then $t\le
a$ or $t\le x\gamma a$ for some $x\in M$ and
$\mu,\gamma\in\Gamma$. Since $M$ is left regular, $a\le
y\mu(a\gamma a)$ for some $y\in M$ and $\mu,\gamma\in\Gamma$. If
$t\le a$, then $t\le a\le y\mu (a\gamma a)$ ($y\in M$,
$\mu,\gamma\in\Gamma)$. If $t\le x\gamma a$, then $t\le x\gamma
a\le x\gamma (y\mu a\gamma a)=(x\gamma y)\mu (a\gamma a)$. In any
case, $t\le z\mu (a\gamma a)$ for some $z\in M$. Hence $t\in
L(a\gamma a)$, and so $L(a)\subseteq L(a\gamma a)$".\\The phrase
"Since $M$ is left regular, $a\le y\mu (a\gamma a)$ for some $y\in
M$ and $\mu,\gamma\in\Gamma$" in it is wrong as the $\gamma$ in
$y\mu (a\gamma a)$ cannot be the same with the $\gamma$ in
$x\gamma a$. The correct is "If $t\in L(a)$, then $t\le a$ or
$t\le x\gamma a$ for some $x\in M$ and $\mu,\gamma\in\Gamma$.
Since $M$ is left regular, $a\le y\mu (a\xi a)$ for some $y\in M$
and $\mu, \xi\in\Gamma$". So the proof of the implication
$(1)\Rightarrow (2)$ in Theorem 1 in [2]
is false, and so is the proof of Theorem 1 given by Kwon and Lee in 
[2].\\
This is a copy from the paper in [1; p. 613, lines 5--11]: "As an
application of the result proved in Theorem 1 in [Kwon, Y. I, Lee,
S. K.: On left regular $po$-$\Gamma$-semigroups. Kangweon-Kyungki
Math. Jour. 6 (1998), No. 2, pp. 149--154] (which is the
References [2] in the present note) we have the following lemmas:

\noindent{\bf Lemma 2.13.} {\it Let $M$ be a $\Gamma$-semigroup.
The following are equivalent:

$(1)$ M is left regular.

$(2)$ $a{\cal L}(a\gamma a)$  for every $a\in M$ and every
$\gamma\in\Gamma$.}\medskip

\noindent{\bf Lemma 2.14.} {\it Let $M$ be a $\Gamma$-semigroup.
The following are equivalent:

$(1)$ M is right regular.

$(2)$ $a{\cal L}(a\gamma a)$  for every $a\in M$ and every
$\gamma\in\Gamma$.}"\\
Although the definition of the $\Gamma$-semigroup considered in
[2] differs from the definition of a $\Gamma$-semigroup considered
in [1], if the Theorem by Kwon-Lee were true, then the two Lemmas
by Hila-Dine would be true as well. It might be also noted that
not only the proof of the Theorem 1 in [2] is wrong, but the
characterization of a $po$-$\Gamma$-semigroup in which $a{\cal
L}(a\gamma a)$ given in [2] is also wrong. The characterization of
a $\Gamma$-semigroup in which $a{\cal L}(a\gamma a)$ (or $a{\cal
R}(a\gamma a))$ given in [1] is wrong as well.

\noindent This note has been submitted in Acta Math. Sinica on 
September 9, 2013 (the date in Greece).
\end{document}